\documentclass[a4paper,11pt]{amsart}
\usepackage{amsfonts,amssymb,eucal,amsmath} 

\newcommand{\cal}[1]{\mathcal #1}
\newcommand{\B}[1]{\mathbb #1}
\newcommand{\Z}{{\B Z}}
\newcommand{\N}{{\B N}}
\newcommand{\s}{{\cal S}} 

\newtheorem{theorem}{Theorem}[section]%
\newtheorem{lemma}[theorem]{Lemma}%
\newtheorem{cor}[theorem]{Corollary}%
\newtheorem{prop}[theorem]{Proposition}%

\newtheorem{rem}[theorem]{Remark}%

\newcommand{\proofend}{\hspace*{\fill}
$\square$ 
 \normalsize\medskip}

\begin{document}

\title[A multivariate arithmetic function]
{A multivariate arithmetic function of combinatorial and topological significance}
\author[V.\,Liskovets]
{\bf Valery A. Liskovets}
\address{V.\,A.\,Liskovets:\quad Institute of Mathematics,
National Academy of Sciences, 220072, Minsk, BELARUS}
\email{liskov@im.bas-net.by}

\date{November 27, 2009} 


\begin{abstract}
We investigate properties of a multivariate function $E(m_1,m_2,\dots,m_r)$,
called {\it orbicyclic}, that arises in enumerative combinatorics in counting
non-isomorphic maps on orientable surfaces. $E(m_1,m_2,\dots,m_r)$ proves
to be multiplicative, and a simple formula for its calculation is provided.
It is shown that the necessary and sufficient conditions for this function
to vanish is equivalent to familiar Harvey's conditions that characterize
possible branching data of finite cyclic automorphism groups of Riemann surfaces.
\end{abstract}

\maketitle

\section{\bf Introduction}

Let $(m_1,m_2,\dots,m_r)$ be a tuple
of $r\geq 0$ positive integers
and $m={\rm lcm}(m_1,m_2,\dots,m_r),$ where $m:=1$ for $r=0$ (an empty
tuple). Introduce the following multivariate function
$$
E=E(m_1,m_2,\dots,m_r):=\frac{1}{M}\sum_{k=1}^M\Phi(k,m_1)\Phi(k,m_2)\cdots\Phi(k,m_r)
\eqno(1)
$$
($E(\emptyset)=1$),
where $m|M,\ M>0,$ and $\Phi(k,n)$ stands for the von Sterneck function:
$$
\Phi(k,n)
:=\frac{\phi(n)}{\phi\bigl(\frac{n}{(k,n)}\bigr)}\,\mu\Bigl(\frac{n}{(k,n)}\Bigr).
\eqno(2)
$$
Here $(k,n)$ denotes the greatest common divisor of $k$ and $n$
and $\mu(n)$ and $\phi(n)$ are the M\"obius and Euler
functions respectively. According to O.\,H\"older (see,
e.g.,~\cite[Ch.\,8]{Ap76},~\cite[Ch.IX]{Si89}),
$\Phi(k,n)$ coincides with the Ramanujan trigonometric sum:
$$
\Phi(k,n)=C_n(k) \eqno(3)
$$
where
$$
C_n(k):=\sum_{d\,({\rm mod}\,n)\atop(d,n)=1}\exp\Bigl(\frac{2\,i\,k\,d}n\Bigr)
$$
with the summation over a reduced residue system modulo $n$.
The function $C_n(k)$ (in the literature it has diverse designations
such as $C(k,n)$) satisfies the familiar Ramanujan's identity (loc. cite):
$$
C_n(k)=\sum_{d|(k,n)}d\,\mu\Bigl(\frac{k}{d}\Bigr).\eqno(4)
$$

Note that $E(m_1,m_2,\dots,m_r)$ does not depend on $M.$ Indeed, by~(2),
$\Phi(k,m_j)$ is a periodic function of $k$ modulo the second
variable and, a fortiori, modulo $m.$ Thus $\prod_j\Phi(k,m_j)$
as a function of $k$ is periodic modulo~$m$ as well. So, in~(1) we
may put $M=m$. Now, $E$ is a symmetric function of its arguments,
and we might speak about the (multi-)set of arguments instead of
a tuple of them. Since $\Phi(k,1)=1$ we may restrict $m_j$ to
values greater $1$, that is,
$$
E(m_1,m_2,\dots,m_{r-1},1)=E(m_1,m_2,\dots,m_{r-1}).\eqno(5)
$$
This plays an important r\^ole in
computational formulae for $E$ and applications.
Tuples of arguments not containing $1$ are called {\it reduced}.

$E(m_1,m_2,\dots,m_r)$ is an essentially multivariate function
in the sense that it is trivial for $r=0,1$. Moreover, for $r=2$\/ it
vanishes for unequal arguments and coincides with the Euler function
otherwise. Later we will see that $E(m_1,m_2,\dots,m_r)$ is always
non-negative and integer.

The function $E$ has been introduced by A.\,Mednykh and
R.\,Nedela in 2004 (see~\cite{MeN04}) in the context of enumerative combinatorics:
it plays a crucial r\^ole in counting maps on orientable surfaces up to
orientation-preserving isomorphism, via the calculation of certain epimorphisms
from the fundamental group of orbifolds onto cyclic groups (see Sect.\,3). Therefore
this ``orbicyclic'' arithmetic function (as we called it in~\cite{Li05})
deserves a detailed investigation by its own right.

Here we study some basic properties of $E(m_1,m_2,\dots,m_r).$
First of all we analyze the prime-power case and establish a simple
sum-free formula for $E(p^{a_1},p^{a_2},\dots,p^{a_r}),$ $p$ prime. It is
determined by three parameters (apart from $p$), denoted $r$ (reduced), $s$
and $v$, only one of which ($s$, the multiplicity of the highest power) is
responsible for its vanishing. In the ge\-ne\-ral case we show that
$E(m_1,m_2,\dots,m_r)$ is a multiplicative function of all its arguments.
Both results provide a simple explicit formula for its calculation.
The most valuable property of $E(m_1,m_2,\dots,m_r)$ established
here is the necessary and sufficient conditions of its
non-vanishing. We show that they are equivalent (due to the
above-mentioned connection with the enumeration of epimorphisms)
to the well-known conditions discovered by W.\,Harvey~\cite{Ha66}
(cf. also~\cite{BuC99}) that specify the possible actions of finite
cyclic automorphism groups of Riemann surfaces. Therefore the
function $E(m_1,m_2,\dots,m_r)$ may be considered as a fruitful
{\sl enumerative refinement} of Harvey's theorem, bringing a new
insight into this theory.

The familiar Jordan arithmetic function also participates in the
mentioned enumeration; we discuss briefly some other links between
this function and enumerative combinatorial group theory.

For the reader's convenience, the paper contains rather numerous (although restricted)
references to relevant publications in the three main to\-pics we deal with:
number theory (arithmetic functions), algebraic to\-po\-lo\-gy (automorphisms of
Riemann surfaces), and algebraic and enumerative combinatorics (map theory
and map enumeration). For a general material concerning these topics and
the notions used in the present paper, the reader is referred to, resp.,~\cite{Ap76},
\cite{BuE90} and~\cite{CoM92}.

\section{\bf Function $E(m_1,m_2,\dots,m_r)$ and its properties}

\subsection{Primary case}
$\Phi(k,n)$ is a multiplicative function of $n$ which is determined
by the following well-known (and easily provable) formula; see, e.g.,~\cite{Mc60}:

\begin{lemma}\label{ram}
For\/ $p$ prime and\/ $a\geq 1,$
$$
\Phi(k,p^a)=\left\{\begin{array}{rl}
(p-1)p^{a-1}& if\quad p^a|k\\
    -p^{a-1}& if\quad p^a\nmid k,\ p^{a-1}|k\\
           0& otherwise.
\end{array}\right.\eqno(6)
$$
\end{lemma}

This is an important result for calculating $E(m_1,m_2,\ldots,m_r)$ explicitly.
As we will see below, the function $E(m_1,m_2,\ldots,m_r)$
is always non-negative, unlike $\Phi(k,n)$ (and it often vanishes like $\Phi(k,n)$).

At first, we consider the prime power case: $m=p^a.$
Let $(m_1,m_2,\ldots,m_r)$ be a reduced tuple.
Then $m_j=p^{a_j},\ j=1,2,\ldots,r,$ where without loss of generality
we assume that
$$
a_1=a_2=\ldots=a_s=a>a_{s+1}\geq a_{s+2}\dots\geq a_r>0,\eqno(7)
$$
where $r\geq s\geq 1.$ Denote
$$
 v:=\sum_{j=2}^r(a_j-1)=\sum_{j=1}^ra_j -r-a+1,\qquad v\geq 0.\eqno(8)
$$
The parameters $s(p)=s$ and $v(p)=v,$ where $p|m,$ can be defined for
an arbitrary tuple of variables $(m_1,m_2,\ldots,m_r)$ as well (see Sect.\,2.3 below).
As we will show, they (together with the corresponding $r(p)$ with respect
to {\sl reduced} tuples) determine the value of $E(m_1,m_2,\ldots,m_r).$

\begin{lemma}\label{prim}
$$
E(p^{a_1},p^{a_2},\ldots,p^{a_r})
=(p-1)^{r-s+1}p^vh_s(p), \eqno(9)
$$
where the multiplicity $s$ and the exponents $a_j$ are subject to~{\rm(7)},
$v$ is determined by~{\rm(8)} and\/ $h_s(x),\ s\geq 1,$ is the following polynomial
of\/ $x$ of degree\/ $s-2$ (for\/ $s>1$):
$$
 h_s(x)=\frac{(x-1)^{s-1}+(-1)^s}{x}. \eqno(10)
$$
\end{lemma}

Proof. We have
$$
E(p^{a_1},\ldots,p^{a_r})=\frac{1}{p^a}\sum_{k=1}^{p^a}\Phi(k,p^{a_1})\cdots\Phi(k,p^{a_r}).
$$
By (6), the first factor in these terms vanishes unless $k=dp^{a-1}.$
So that we may restrict ourselves to such $k$ only, where $d=1,2,\dots,p.$ Again by (6),
we get (for all $s$ including $s=1$)
$$
E(p^{a_1},\ldots,p^{a_r})
=\frac{1}{p^a}\Big(\sum_{d=1}^{p-1}\Phi(dp^{a-1},p^{a_1})\cdots\Phi(dp^{a-1},p^{a_r})
$$
$$
+\,\Phi(p^a,p^{a_1})\cdots\Phi(p^a,p^{a_r})\Big)
$$
$$
=\frac{1}{p^a}\!\Big((p-1)\Phi(p^{a-1},p^{a_1})\cdots\Phi(p^{a-1},p^{a_r})
+\Phi(p^a,p^{a_1})\cdots\Phi(p^a,p^{a_r})\Big)
$$
$$
=\frac{1}{p^a}\Big((p-1)(-p^{a-1})^s(p-1)^{r-s}p^{a_{s+1}-1}\cdots p^{a_r-1}
+(p-1)^rp^{a_1-1}\cdots p^{a_r-1}\Big)
$$
$$
=\frac{1}{p^a}\Big((p-1)(-1)^sp^{(a-1)s}p^{\sum_{j=s+1}^r(a_j-1)}(p-1)^{r-s}+(p-1)^rp^{\sum_{j=1}^r(a_j-1)}\Big)
$$
$$
=\frac{1}{p^a}\Big((p-1)(-1)^sp^{\sum_{j=1}^r(a_j-1)}(p-1)^{r-s}+(p-1)^rp^{\sum_{j=1}^r(a_j-1)}\Big)
$$
$$
=\frac{1}{p^a}p^{\sum_{j=1}^r(a_j-1)}(p-1)^{r-s+1}\Big((-1)^s+(p-1)^{s-1}\Big)
$$
$$
=p^{-1+\sum_{j=2}^r(a_j-1)}(p-1)^{r-s+1}\Big((-1)^s+(p-1)^{s-1}\Big)
$$
$$
=(p-1)^{r-s+1}p^v\bigg(\frac{(-1)^s+(p-1)^{s-1}}{p}\bigg).
$$
\proofend

In particular,
$$
\begin{array}{ll}
&h_1(x)=0,\\
&h_2(x)=1,\\
&h_3(x)=x-2,\\
&h_4(x)=x^2-3x+3,\\
&h_5(x)=(x-2)(x^2-2x+2).
\end{array}
$$
Note that by~(10), $(x-2)|h_s(x)$ for odd $s$ and $h_s(2)=1$ for even $s.$
Clearly $h_s(p)\geq 0$ for any $s$ and $p\geq 2.$ Moreover, $h_s(p)=0$ if and only
if $s=1$ or $p=2$ and $s$ is odd.

\begin{cor}\label{van1}
{\bf (1)} $E(p^{a_1},p^{a_2},\ldots,p^{a_r})$ is non-negative and integer.

{\bf (2)} $E(p^{a_1},p^{a_2},\ldots,p^{a_r})$ vanishes if and only if\/
$s=1$ or\/ $p=2$ and\/ $s$ is odd.

{\bf (3)} Let $p$ be odd. $p\nmid E(p^{a_1},p^{a_2},\ldots,p^{a_r})$ if and only if
$a=1$ and $r>1;$ in this case, $s=r.$

{\bf (4)} $E(p^{a_1},p^{a_2},\ldots,p^{a_r})$ is odd if and only if $p=2,a=1$
and $s=r$ is even.
\end{cor}

Proof. Claims~(1) and~(2) are immediate from formulae~(9) and~(10).

Claims~(3). Suppose that $p$ does not divide $E(p^{a_1},p^{a_2},\ldots,p^{a_r}).$
Then $E(p^{a_1},p^{a_2},\ldots,p^{a_r})\neq 0$ and by formula~(9), $v=0.$
It follows from~(8) that $a_2=\dots=a_r=1.$ But $a_1=1$ as well and $r>1$,
otherwise $s=1$ and then by Claim~(2), $E(p^{a_1},p^{a_2},\ldots,p^{a_r})=0.$
Thus, $a=a_1=1$ and $s=r.$

On the contrary, if $a=1$ then $a_1=a_2=\dots=a_r=1$ as well. Thus, $s=r$ and by~(8),
$v=0.$ But in the right-hand side of formula~(9), neither $p-1$ nor $h_s(p)$ is divisible
by $p$ for $p>2$ and $s>1.$ Hence if $r>1$ then $p\nmid E(p^{a_1},p^{a_2},\ldots,p^{a_r}).$

Claims~(4). $E(p^{a_1},p^{a_2},\ldots,p^{a_r})$ is even for odd $p$ since always $r-s+1>0$ and,
thus, the factor $(p-1)^{r-s+1}$ in the right-hand side of formula~(9) is even. Now
suppose $p=2.$ We have the same situation as in the proof of Claim~(3) with the only distinction
that $2|h_s(2)$ iff $s$ is odd.
\proofend

\subsection{Multiplicativity}
Now we turn to the general case. Take a tuple $(m_1,m_2,\ldots,m_r)$
and let $a_j(p)\geq 0$ denote the exponent with which the prime $p$ divides
$m_j$. For a prime $p|m$, denote
$$
 \{m_1,m_2,\ldots,m_r\}_p=(p^{a_1(p)},p^{a_2(p)},\ldots,p^{a_r(p)}).
 \eqno(11)
$$
Now denote
$$
E_p(m_1,m_2,\ldots,m_r):=E(\{m_1,m_2,\ldots,m_r\}_p)=E(p^{a_1(p)},p^{a_2(p)},\ldots,p^{a_r(p)}).
\eqno(12)
$$
Then
$$
 E_p(m_1,m_2,\ldots,m_r)=E(\langle p^{a_1(p)},p^{a_2(p)},\ldots,p^{a_r(p)}\rangle),
 \eqno(13)
$$
where
the angled brackets mean the removal of all the arguments equal to $1$,
that is, the right-hand function contains a reduced set of arguments.

\begin{prop}\label{mult1}
$E(\emptyset)=E(1,1,\ldots,1)=1,$ and for\/ $m>1$,
$$
 E(m_1,m_2,\ldots,m_r)=\prod_{p|m\ \rm prime}E_p(m_1,m_2,\ldots,m_r).
\eqno(14)
$$
In other words, $E(m_1,m_2,\ldots,m_r)$ is a multiplicative
function. 
\end{prop}

An arithmetic function $g=g(m_1,m_2,\dots,m_r)$ of $r\geq 1$ arguments
is called {\it semi-multiplicative} (see, e.g.,~\cite{Ha88,Re66})
if
$$
g(m_1,m_2,\dots,m_r)=g(m_1',m_2',\dots,m_r')\cdot g(m_1'',m_2'',\dots,m_r'')
$$
whenever $m_j=m_j'm_j''$, $j=1,2,\dots,r,$ and $(M',M'')=1$
where $M'=\prod_jm_j'$ and $M''=\prod_jm_j''$. This function is called
{\it multiplicative} if, moreover, $g(\underbrace{1,1,\dots,1}_r)=1$.
Two most known examples of symmetric multivariate multiplicative functions are
lcm$()$ and gcd$()$.

\medskip
$\Phi(k,n)$ is multiplicative as a function of $n$ and is periodic in $k$
modulo $n$ (it possesses, in fact, stronger properties but we do not need
to use them here). Therefore
Proposition~\ref{mult1} is straightforward from the following more
general assertion.

\begin{lemma}\label{periodic}
Let a bivariate arithmetic function $f(k,n)$ be semi-multiplicative in
the argument\/ $n$ and periodic in\/ $k$ modulo\/ $n.$ Given natural
numbers $M,m_1,m_2,\ldots,m_r,$ where $m_j|M$ for all $j$, the function
$$
F(m_1,m_2,\dots,m_r)=\frac{1}{M}\sum_{k=1}^M f(k,m_1)f(k,m_2)\cdots f(k,m_r)\eqno(15)
$$
is semi-multiplicative (with respect to all its arguments).
\end{lemma}

Proof.
By the same reasons as in the Introduction, $F(m_1,m_2,\dots,m_r)$
does not depend on $M$ (provided $m_j|M,\ j=1,2,\dots,r$). So that,
without loss of generality, we set $M=\prod_j^rm_j.$ Suppose
$m_j=m_j'm_j''$, $j=1,2,\dots,r,$ where $(m_i',m_j'')=1$ for all $i,j.$
Denoting $M'=\prod_j^rm_j'$ and  $M''=\prod_j^rm_j''$ consider the product
$$
F(m_1',\dots,m_r')\cdot F(m_1'',\dots,m_r'')
$$
$$
=\frac{1}{M'}\sum_{k'=1}^{M'}f(k',m_1')\cdots f(k',m_r')
\cdot\frac{1}{M''}\sum_{k''=1}^{M''}f(k'',m_1'')\cdots f(k'',m_r'').
$$
We claim that, regardless of the choice of all $m_j=m_j'm_j''$, the equality
$F(m_1,\dots,m_r)=F(m_1',\dots,m_r')\cdot F(m_1'',\dots,m_r'')$ holds.
Notice first that like $F(m_1,\dots,m_r)$, the product
$F(m_1',\dots,m_r')\cdot F(m_1'',\dots,m_r'')$ contains $M=M'M''$ terms.
Thus, all we need is to establish a bijection between the terms of the
two functions and to show the equality of the corresponding terms.

Since the function $f(k,m)$ is periodic in\/ $k$ modulo\/ $m,$ it is periodic modulo
$M,$ as well, if $m|M,$ i.e., $f(k_1,m)=f(k_2,m)$ whenever $k_1\equiv k_2\ ({\rm mod}\ M)$.

Consider a 'generic' term in $F(m_1',\dots,m_r')\cdot F(m_1'',\dots,m_r''):$
$$
t(k',k''):=f(k',m_1')\cdots f(k',m_r')\cdot f(k'',m_1'')\cdots f(k'',m_r'').
$$
We look for $k$ such that $t(k',k'')=t(k,k).$ Take $k$ satisfying
$$
\begin{array}{ll}
&k\equiv k'\ ({\rm mod}\ M')\\
&k\equiv k''\ ({\rm mod}\ M'').
\end{array}
$$
Since $(M',M'')=1,$ by the Chinese remainder theorem, such a $k$
does exist and is unique modulo $M.$ Now, due to the periodic property,
$f(k',m_j')=f(k,m_j')$ and $f(k'',m_j'')=f(k,m_j'')$ for $j=1,2,\dots,r.$ Hence
$t(k',k'')=t(k,k)$ as required. Therefore by the semi-multiplicativity of $f(k,m)$
in $m$ we have
$$
t(k',k'')=f(k,m_1')\cdots f(k,m_r')\cdot f(k,m_1'')\cdots f(k,m_r'')=
f(k,m_1)\cdots f(k,m_r).
$$

It is clear that the established correspondence between $k$ and pairs $k',k''$
is a bijection between the sets $[1,M]$ and $[1,M']\times[1,M''],$
what gives rise to the required bijection between the terms of
$F(m_1,m_2,\dots,m_r)$ and those of $F(m_1',\dots,m_r')\cdot F(m_1'',\dots,m_r'').$
\proofend

\begin{cor}\label{split}
If $m_1=m'_1m''_1,$ where
$(m'_1,m''_1)=1,$
then
$$
E(m_1,m_2,\ldots,m_r)=E(m'_1,m''_1,m_2,\ldots,m_r).
$$
\end{cor}

Indeed,
$E_p(m_1,m_2,\ldots,m_r)=
{E}_p(\langle m'_1,m''_1,m_2,\ldots,m_r\rangle)$
for all prime $p.$
\proofend

Therefore, increasing $r$, one can split the arguments of $E(m_1,m_2,\ldots,m_r)$
into their primary factors:
$$
E(m_1,m_2,\ldots,m_r)
=E\bigl(p^{a_j(p)}:\,\, p|m\,\, {\rm prime},\,\, j=1,2,\ldots,r,\, a_j(p)\geq 1\bigr).
$$

\begin{cor}\label{posit}
The values of $E(m_1,m_2,\ldots,m_r)$ are non-negative integers. \proofend
\end{cor}

\subsection{Main formulae}
Given a tuple $(m_1,m_2,\ldots,m_r)$ with
$$
{\rm lcm}(m_1,m_2,\ldots,m_r)=m=\prod\limits_{p\ \rm prime}p^{a(p)},
\eqno(16)
$$
define for $p|m,$ the parameters $s(p)$
and $v(p)$ that generalize the ones introduced in formulae~(7) and~(8):
$$
 s(p):=|\{j:\ a_j(p)=a(p),\ j=1,2,\ldots,r\}| \eqno(17)
$$
and
$$
 v(p):=\sum_{j=1,2,\ldots,r\atop a_j(p)\geq 1}(a_j(p)-1)-a+1.\eqno(18)
$$
Moreover, we count the arguments $m_j$ divisible by $p,$ that is,
the ones with $a_j(p)>0:$
$$
 r(p):=|\{m_j:\ p|m_j,\ j=1,2,\ldots,r\}|.\eqno(19)
$$

The next theorem follows directly from~(9) and~(14) and gives rise to
an explicit alternating-free formula for calculating $E(m_1,m_2,\ldots,m_r).$

\begin{theorem}\label{main}
$$
 E_p(m_1,m_2,\ldots,m_r)
 =(p-1)^{r(p)-s(p)+1}p^{v(p)}h_{s(p)}(p) \eqno(20)
$$
and
$$
 E(m_1,m_2,\ldots,m_r)
=\prod_{p|m\ \rm prime}(p-1)^{r(p)-s(p)+1}p^{v(p)}h_{s(p)}(p),\eqno(21)
$$
where the parameters\/ $s(p),v(p)$ and\/ $r(p)$ are defined, respectively,
by formulae~{\rm(17)},~{\rm(18)} and~{\rm(19)} and the polynomial\/
$h_s(x)$ is defined by~{\rm(10)}. \proofend
\end{theorem}

\subsection{Further properties}

According to~(21), the value of the function $E(m_1,m_2,\ldots,m_r)$
is determined by the set of prime divisors $p|m$ and the parameters
$s(p),v(p)$ and $r(p),$ where $r\geq r(p)\geq s(p)\geq 1$ and $v(p)\geq 0.$
In particular, $E(m_1,m_2,\ldots,m_r)$ does not depend directly on $a(p)$
and $m$ (indirectly, however, $a(p)$ contributes into $v(p)$).
Note also that $E(m_1,m_2,\ldots,m_r)$ does not depend on $r(2)$ as
formulae~(20) and~(18) show; so that if $4|m,$ we may ignore the contributors
$2^1$, that is, $a_j(2)=1.$

\begin{cor}\label{divis}
{\bf (1)} $\phi(m)$\/ divides\/ $E(m_1,m_2,\ldots,m_r).$

{\bf (2)} $E(m_1,m_2,\ldots,m_r)=\phi(m)$ if and only if for every prime\/
$p|m,$ one of the following conditions holds:
$$
\langle\{m_1,m_2,\ldots,m_r\}_p\rangle=(p^{a(p)},p^{a(p)})
$$
or $p=3$ and
$$
 \langle\{m_1,m_2,\ldots,m_r\}_3\rangle=(3,3,3)
$$
or $p=2$ and (up to reordering)
$$
 \langle\{m_1,m_2,\ldots,m_r\}_2\rangle
=(2^{a(2)},2^{a(2)},\underbrace{2,2,\ldots,2}_{r(2)-2}),
$$
where\/ $a(2)\geq 1,\,r(2)\geq 3$ and $r(2)$ is even if $a(2)=1.$
\end{cor}

Proof. Claim~(1). Recall that
$\phi(p^{a(p)})=(p-1)p^{a(p)-1}$. Now, for any $p|m,$
it follows directly from the definitions that $r(p)-s(p)+1\geq 1.$ Besides,
$v(p)\geq a(p)-1$ if $E_p(m_1,m_2,\ldots,m_r)\neq 0$ since in this case,
$s(p)>1$ by Corr.~\ref{van1}~(2). Thus, in~(18), the term corresponding to $j=2$
is equal to $a(p)-1.$
Therefore, by~(20), $\phi(p^{a(p)})|E_p(m_1,m_2,\ldots,m_r),$
and we are done by~(14) and the multiplicativity of $\phi.$

Claim~(2). It is clear from~(10) that $p$ does not divide $h_s(p)$ for $s>1,$
nor does $p-1$ for $p>2.$ Therefore by~(21), $E(m_1,m_2,\ldots,m_r)=\phi(m)$
if and only if
$$
 r(p)-s(p)+1=1,
$$
$$
 v(p)=a(p)-1
$$
and
$$
 h_{s(p)}(p)=1
$$
for every prime $p|m.$ Suppose $p\geq 3.$ The first equality implies that
$\langle\{m_1,m_2,\ldots,m_r\}_p\rangle=(p^{a(p)},p^{a(p)},\ldots,p^{a(p)})$
whence $v(p)=(r(p)-1)(a(p)-1).$ Then the second equality implies that
$r(p)=2$ or $a(p)=1.$ If $r(p)=s(p)=2$ then we get $h_{s(p)}(x)=1.$
If $r(p)\neq 2$ but $a(p)=1,$ then $h_s(p)=1$ is possible only in one
exceptional case: $h_3(3)=1.$ Indeed, $h_s(x)=1$ can be represented
as the equation $y^{s-1}-y-1+(-1)^s=0,$ where $y=x-1.$ It
has only one integer solution greater 1: $y=2$ and $s=3.$
This solution corresponds to the tuple $(3,3,3)$ as claimed.

Finally, consider $p=2$ if $2|m.\quad h_s(2)=1$ if and only if $s$ is even.
Now $v(2)=a(2)-1$, what means that
$\langle\{m_1,m_2,\ldots,m_r\}_2\rangle=(2^{a(2)},2^{a(2)},2,2,\ldots,2)$
with $a(2)\geq 1.$ If $a(2)>1,$ then $s(2)=2$ is even; if $a(2)=1,$ then
$s(2)=r(2)$ has to be even by Corr.~\ref{van1}~(4).
\proofend

Let us consider in more detail the behavior of the function $E$ for $r\leq 3$
arguments.

\begin{cor}\par
{\bf (1)} $E(m)\neq 0$ if and only if\/ $m=1.$
\label{1-2-3}

{\bf (2)} $E(m_1,m_2)\neq 0$ if and only if\/ $m_1=m_2.$

{\bf (3)} $E(m_1,m_2,m_3)\neq 0$ if and only if
for any prime\/ $p|m,$ the numbers $m_1,m_2$ and $m_3$ (in some order) are divided
by\/ $p$ with the exponents\/ $a(p)\geq b(p)\geq c(p)$ that meet one of the following
conditions:
\begin{itemize}
\item $a(p)=b(p)>c(p)>0$;
\item $a(p)=b(p)>c(p)=0$;
\item $a(p)=b(p)=c(p)>0$\ and $p>2$.
\end{itemize}
\end{cor}

Proof. The first claim is obvious. As to the second one,
if $m_1\neq m_2$ then there is a prime $p|m$ such that $s(p)=1.$
Therefore the last factor in formula~(21) vanishes.

In claim~(3), if ${a(p)>b(p)\geq c(p)}$ for some $p|m,$
then $s(p)=1$ and the last factor $h_s$ in formula~(20) vanishes.
It also vanishes if $2|m,p=2$ and $s(2)=3,$ that is, $a(2)=b(2)=c(2).$
\proofend

Now, $E_p$ is evaluated by formula~(9), and we obtain for the last case
of Corr.~\ref{1-2-3},
$$
E_p(m_1,m_2,m_3)=\left\{\begin{array}{ll}
(p-1)^2p^{a(p)+c(p)-2}& {\rm if}\ a(p)=b(p)>c(p)>0\\
(p-1)p^{a(p)-1}       & {\rm if}\ a(p)=b(p)>c(p)=0\\
(p-1)(p-2)p^{2a(p)-2} & {\rm if}\ a(p)=b(p)=c(p)>0
\end{array}\right\}. \eqno(22)
$$

Given $m=\prod_{p|m}p^{a(p)},$ Corr.\,\ref{1-2-3}~(3) and
Prop.\,\ref{mult1} make it possible to easily observe all triples
$m_1,m_2,m_3$ for which $E(m_1,m_2,m_3)\neq 0.$ Namely, for any
prime $p|m$ we first form a triple $T_p=\{p^{a(p)},p^{a(p)},p^{c(p)}\}$
with an arbitrary integer $c(p),\ 0\leq c(p)\leq a(p).$ Then we form
three sets $\cal M_1,\cal M_2$ and $\cal M_3$ by distributing the elements of
each $T_p$ arbitrarily by them (so that, there are three different ways
if $c(p)<a(p)$, and only one way otherwise).
Now we obtain the desired numbers $m_j$ as the products of the elements in
the corresponding $\cal M_j$:\,
$m_j:=\prod_{p|m}p^b,\, j=1,2,3,$ where $p^b\in\cal M_j$.
It follows that there are totally $\prod_{p|m}(3a(p)+1)$ such
ordered triples.

\begin{cor}
$f_r(m):=E(\underbrace{m,m,\dots,m}_r)$
is the multiplicative function of $m$ determined by the formula
\label{mmm}
$$
 f_r(p^a)=(p-1)p^{(r-1)(a-1)}h_r(p),\quad p\ \ {\rm prime},\ a>0. \eqno(23)
$$
\end{cor}
\proofend

According to~(23), $f_r(m),\, m>1,$ can be represented as follows:
$$
 f_r(m)=m^{r-1}\prod_{p|m\ \rm prime}\frac{(p-1)h_r(p)}{p^{r-1}}.
$$

It follows that $E(\underbrace{m,m,\dots,m}_r)=0$
if and only if $m$ is even and $r$ is odd or $r=1$ and $m>1$.
Besides~\cite{Ni53},
$$
 E(m,m)=\phi(m) \eqno(24)
$$
(in other words, $f_2(m)=\phi(m)$). This bivariate instance of $E$ is the
only non-trivial particular specimen of the function defined by formula~(1)
that I managed to find in the number-theoretic literature.

Likewise, relying upon formula~(21), we could investigate other general properties
of the function $E(m_1,m_2,\ldots,m_r).$
One such property significant in topological applications, namely non-vanishing,
will be considered below in Theorem~\ref{van2} (in the particular case $m=p^a$ we have
already encountered it in Corr.~\ref{van1}~(2) and later).

\section{\bf Combinatorial and topological motivations and applications}

\subsection{Linear congruences (restricted partitions)}
The following statement provides a simple combinatorial interpretation
for the function $E:$

\begin{lemma} {\rm\cite{MeN04}}\label{congr}
Let $M$ be a natural number and\/ $m_1,m_2,\ldots,m_r$ divisors of $M.$
Denote
${d_1=\frac{M}{m_1},d_2=\frac{M}{m_2},\ldots,d_r=\frac{M}{m_r}}.$
Then the number of solutions of the system of equations
$$
\left.
\begin{array}{rl}
x_1+x_2+\cdots+x_r\equiv &0\ ({\rm mod}\,M)\\
(x_1,\,M)=&d_1\\
\ldots&\\
(x_r,\,M)=&d_r
\end{array}\right\} \eqno(25)
$$
in integers modulo~$M$ does not depend on $M$ and is equal to
$E=E(m_1,m_2,\ldots,m_r).$
\end{lemma}

One could use Lemma~\ref{congr} to give another proof of Proposition~\ref{mult1}.

\subsection{Jordan's function}
A more profound combinatorial interpretation for the function
$E(m_1,m_2,\dots,m_r)$ is presented below in Theorem~\ref{epi}
and combines $E$ with the classical Jordan function.
Recall~\cite{Ap76,Co59,Sc99} that the Jordan function of order $k$
is defined by
$$
 \phi_k(n)=\sum_{d|n}d^k\mu\Bigl(\frac{n}{d}\Bigr) \eqno(26)
$$
or, equivalently, as a multiplicative function, by
$$
 \phi_k(n)=n^k\prod_{p|n\, \rm prime}(1-p^{-k}).\eqno(27)
$$
In particular, $\phi_1=\phi$ (Euler's totient). Moreover,
$\phi_0(1)=1$ and $\phi_0(n)=0$ for $n>1.$ It follows that
$$
 \phi(n)|\ \phi_k(n)
$$
for $k\geq 1.$

As can be noticed in the literature on arithmetic functions,
it is rather typical to research the Ramanujan sums and their
generalizations jointly with the Jordan functions (cf.,
e.g.,~\cite{Co59,Mc60}).

\begin{rem}\label{rem1}
{\rm It is worth observing
that the Jordan function participates (expressly or implicitly) in reductive enumeration
formulae for conjugacy classes of subgroups
of some finitely ge\-ne\-ra\-ted groups. The oldest formula of this
kind was derived by the author~\cite{Li71} and can be represented as follows:
$$
N_{F_r}(n)
=\frac{1}{n}\sum\limits_{d|n}\phi_{(r-1)d+1}\Bigl(\frac{n}{d}\Bigr)M_{F_r}(d),
\eqno(28)
$$
where $F_r$ is a free group of rank $r$, $M_G(n)$ denotes the number of
subgroups of index $n$ in a group $G$ and $N_G(n)$ denotes the number
of conjugacy classes of such subgroups. Subsequent similar formulae obtained
by A.\,Mednykh in the 1980s and later for $n$-index subgroups of
the fundamental groups of closed surfaces (see~\cite{Me06,Me08,KMN0x}
and~\cite[Sect.\,2 and 5]{Li98}) also contain
$\phi_k(n/d)$ as a factor, where $k$ is a linear function of $d.$
A considerable r\^ole of the Jordan function in this context has been
realized only recently.}
\end{rem}

\subsection{Orbifolds and cyclic automorphism groups of Riemann surfaces}
In order to describe more sophisticated applications of $E$ we need to remind
some notions of algebraic topology.
Generally for automorphisms of Riemann surfaces and orbifolds we refer
to~\cite[Sect.\,2]{Sc83} and~\cite[Ch.\,3]{BuE90}. By an {\it orbifold} we mean
here the quotient space of a closed orientable surface with respect to an action
of a finite group of orientation-preserving automorphisms.
Denote by ${\rm Orb}(\s_\gamma/G)$ the set of orbifolds arising as the
quotient spaces by the actions of the group $G$ on a Riemann surface $\s_\gamma$
of genus $\gamma$.
Any orbifold $\Omega\in {\rm Orb}(\s_\gamma/G)$ is a closed
surface $\s_g$ with a finite set of distinguished branch points; a surface
is a particular case of an orbifold, with the empty set of branch points.
$\Omega$ is characterized by its {\it signature} $(g;m_1,m_2,\dots,m_r),\ r\geq 0,$
where ${2\leq m_j\leq\ell},\ j=1,2,\ldots,r,$ are the orders of its branch points
and $|G|=\ell.$ We denote $\Omega=\Omega(g;m_1,m_2,\dots,m_r)$. Following Mednykh
and Nedela~\cite{MeN04}, by {\it cyclic orbifolds} we mean orbifolds
corresponding to the cyclic groups $G=\Z_\ell.$

Given an orbifold $\Omega=\Omega(g;m_1,m_2,\ldots,m_r),$
define the orbifold fundamental group $\pi_1(\Omega)$ to be the (Fuchsian) group
generated by $2g+r$ gene\-ra\-tors $x_1,y_1,x_2,y_2,\dots,x_g,y_g$ and
$z_1,z_2,\dots,z_r$ and satisfying the relations
$$
\prod_{i=1}^g[x_i,y_i]\prod_{j=1}^r z_j=1\quad {\rm and}\quad
z_j^{m_j}=1,\quad j=1,\ldots,r, \eqno(29)
$$
where $[x,y]=xyx^{-1}y^{-1}$.

Actions of a group $G$ on a surface naturally correspond to epimorphisms
from the fundamental group of the corresponding orbifold onto $G$.

\subsection{Epimorphisms $\pi_1(\Omega)\to\Z_\ell$}
An epimorphism from $\pi_1(\Omega)$ onto a cyclic group of order $\ell$
is called {\it order-preserving} if it preserves the orders of the
periodical generators $z_j$, $j=1,\ldots,r$. Equivalently, order-preserving
epimorphisms have torsion-free kernels; in the literature,
such epimorphisms are often called {\it smooth}, see, e.g.,~\cite{MM98}.
We denote by ${\rm Epi}_{\rm o}(\pi_1(\Omega),\Z_\ell)$ the set of
order-preserving epimorphisms $\pi_1(\Omega)\to\Z_\ell$. Lemma~\ref{congr}
makes it possible to find their number $|{\rm Epi}_{\rm o}(\pi_1(\Omega),\Z_\ell)|:$

\begin{theorem} {\rm\cite{MeN04}}\label{epi}
The number of order-preserving epimorphisms
from the fundamental group\/ $\pi_1(\Omega)$ of the cyclic orbifold\/
$\Omega=\Omega(g;m_1,m_2,\ldots,m_r)\in{\rm Orb}(\s_\gamma/\Z_\ell)$ onto the cyclic
group\/ $\Z_\ell$ is expressed by the following formula:
$$
 |{\rm Epi}_{\rm o}(\pi_1(\Omega),\Z_\ell)|
 =m^{2g}\phi_{2g}(\ell/m)\cdot E(m_1,m_2,\dots,m_r),\eqno(30)
$$
where\/ $m={\rm lcm}(m_1,m_2,\dots,m_r)$ and\/ $\phi_{2g}(m)$ is the Jordan
function of order\/ $2g$.
\end{theorem}

In particular, for $r=0,$
$$
|{\rm Epi}_{\rm o}(\pi_1(\Omega(g;\emptyset)),\Z_\ell)|=\phi_{2g}(\ell)
$$
and for $g=0,$
$$
|{\rm Epi}_{\rm o}(\pi_1(\Omega(0;m_1,m_2,\dots,m_r)),\Z_\ell)|=E(m_1,m_2,\dots,m_r)
$$
if, and only if, $m=\ell.$

Here and subsequently we follow the {\bf convention} that an arithmetic function
vanishes for non-integer arguments
(besides, we pay the reader's attention to a minor change
in designations with respect to $E(m_1,\dots,m_r)$: by some reasons, from now on,
the letter $\ell$ is used instead of $M$).

Due to formula~(30), we call $E(m_1,m_2,\dots,m_r)$ the
{\it orbicyclic} (multivariate arithmetic) function\footnote{Curiously,
a similar but different term ({\sl orbicycle index polynomial}) has been
introduced in~\cite{BlD05}, again in connection with orbit enumeration.}.

According to Corollary~\ref{divis} (1) and formula~(30),
$|{\rm Epi}_{\rm o}(\pi_1(\Omega),\Z_\ell)|$ is divisible by $\phi(\ell),$
what is obvious combinatorially since $\phi(\ell)$ is the number of
primitive elements (units) in the group $\Z_\ell$. Moreover, by~(30),
$$
 |{\rm Epi}_{\rm o}(\pi_1(\Omega),\Z_\ell)|=\phi(\ell)
$$
if and only if $g=0,\,m=\ell$ and $m_1,m_2,\dots,m_r$ satisfy
the conditions of Corollary~\ref{divis}~(2).
Of course in some cases, this can be established directly from~(30).
For example, if $g=0$ and $r=2,$ then it is clear from~(30) that
$\Omega=\Omega(0,m_1,m_2)$ exist if and only if $m_1=m_2=m,$ and in this case
$\pi_1(\Omega)=\Z_m.$ Now, order-preserving epimorphisms from
$\Z_m$ onto $\Z_\ell$ exist if and only if $m=\ell.$ Finally it is
obvious that $|{\rm Epi}_{\rm o}(\Z_m,\Z_m)|=\phi(m).$

\subsection{Riemann--Hurwitz equation}
The signature of any orbifold
$\Omega(g;m_1,m_2,\dots,m_r)\in{\rm Orb}(\s_\gamma/G)$
sa\-ti\-s\-fies the famous Riemann--Hurwitz equation:
$$
2-2\gamma=\ell\Bigl(2-2g-\sum_{j=1}^r\Bigl(1-\frac1{m_j}\Bigr)\Bigr), \eqno({\rm RH})
$$
where $\ell=|G|.$ It is clear
that in~(RH),
$$
 g\leq \gamma. \eqno(31)
$$

The following well-known inequalities can be easily deduced
{\sl directly} from the Riemann--Hurwitz equation
(cf.~\cite[4.14.27]{ZVC80}):

\begin{prop}\label{wim}
{\bf (1)} If\/ $\Omega(g;m_1,m_2,\dots,m_r)\in{\rm Orb}(\s_\gamma/G),\ |G|=\ell,$
then
$$
 \ell\leq 4\gamma+2\quad{\rm for}\quad g\geq 2. \eqno(32)
$$
Thus (see~{\rm\cite{Ha66,Br90}}),
$$
\ell\leq 4\gamma+2\quad{\rm for}\quad\gamma\geq 2\quad{\rm(Wiman,\,1895)}.\eqno(33)
$$

{\bf (2)} For all $\gamma\geq 0,$
$$
r\leq 2\gamma+2. \eqno(34)
$$
\end{prop}

Proof. Denoting $f=\sum_{j=1}^r(1-1/m_j)$,~(RH) can be rewritten as follows:
$4\gamma+2={2\ell(2g+f-2)+6}.$ If\/ $2\ell(2g+f-2)+6<\ell$, then
$\ell(4g+2f-5)+6<0$ and ${4g+2f-5<0}$,
which is possible only for $g=0,1$ since $f\geq 0.$

Wiman's inequality~(33) is immediate from~(31) and~(32).

Given $\gamma$, it is easy to see by~(RH) that $r$ is maximal
whenever $g=0$ and $m_j=m=\ell=2$ for all $j,$ in which case
$r=2\gamma+2.$
\proofend

\begin{rem}\label{rem4}
{\rm {\bf (1)} There is a subtler restriction on $\ell$~\cite{Br90}:
$$
 \phi(\ell)\leq 2\gamma\quad{\rm for}\quad\gamma\geq 2. \eqno(35)
$$
For example, it excludes the order $\ell=9$ for $\gamma=2$,
which satisfy~(33). The bound~(35) excludes also prime $\ell>2\gamma+1.$

\medskip
{\bf (2)} Wiman's bound~{\rm (33)} is also valid for $\gamma=1$ provided that
$r\geq 1.$}
\end{rem}

\subsection{Harvey's conditions for automorphisms of Riemann surfaces}

\begin{theorem} {\rm\cite{Ha66}}\label{harv}
There exists an orbifold\/
$\Omega(g;m_1,m_2,\dots,m_r)\in {\rm Orb}(\s_\gamma/\Z_\ell),$ where
$\ell\geq 2$ and $\gamma\geq 2,$
if and only if its parameters satisfy the Riemann--Hurwitz equation\ {\rm (RH)}
and the following conditions:

\begin{itemize}
\item[(H1)] ${\rm lcm}(m_1,m_2,\dots,m_{j-1},m_{j+1},\dots,m_r)=m$\quad for every
$j=1,2,\dots,r,$ where $m={\rm lcm}(m_1,m_2,\dots,m_r)$\quad (the
{\rm lcm}-condi\-ti\-on);
\item[(H2)] $m$ divides $\ell$, and $m=\ell$ if $g=0$;
\item[(H3)] $r\neq 1$, and $r\geq 3$ if $g=0$;
\item[(H4)] if\/ $m$ is even, then the number of $m_j$
divisible by the maximal power of $2$ dividing $m$ is even.
\end{itemize}
\end{theorem}

\begin{rem}\label{rem3}
{\rm Theorem~{\rm\ref{harv}} is valid for $\gamma=0,1$ as well with
the following condition that supplements~{\rm (H3)}~\cite{MeN04}:
\medskip

\begin{itemize}
\item[(H3a)] $r=2$ if $\gamma=0$, and $r\in\{0,3,4\}$ if $\gamma=1$.
\end{itemize}}
\end{rem}

\subsection{Further applications of the function $E$}

Note that $f_r(m)$ (formula~(23)) gives rise to the number of solutions
of the congruences~(25) when $d_1=d_2=\ldots=d_r=1.$ This is a
particular case of the system consi\-de\-red and solved in~\cite{NiV54} and later
in~\cite{Co55} (see also~\cite{Co59}), where an arbitrary $n$ stands in place of $0$
in the right hand side of the congruence.
This was done in terms of the Ramanujan sums as well.
A related enumeration problem was considered in~\cite{Ug93}.
\medskip

A simple enumerative proof can be obtained for another familiar result
of Harvey that supplements the bound~(33):

\begin{cor}{\rm\cite{Ha66}}\label{harv1}
${\rm Orb}(\s_\gamma/\Z_{4\gamma+2})\neq\emptyset$
for every $\gamma.$ In other words, an orientable surface of genus $\gamma$
possesses an (orientation-preserving) automorphism of order $\ell=4\gamma+2.$
\end{cor}

Proof. Indeed,
$\Omega={\Omega(0;4\gamma+2,2\gamma+1,2)\in{\rm Orb}(\s_\gamma/\Z_{4\gamma+2})}$
exists since the triple $m_1={4\gamma+2},m_2={2\gamma+1},m_3=2$ satisfies
the Riemann--Hurwitz equation~(RH) with $\ell=m={4\gamma+2}$ and $g=0.$
Now by~(30),
$$
|{\rm Epi}_{\rm o}(\pi_1(\Omega),\Z_{4\gamma+2})|=E(4\gamma+2,2\gamma+1,2)
=E(4\gamma+2,4\gamma+2)=\phi(4\gamma+2)>0.
$$ \proofend

Note that $\Omega=\Omega(g;m_1,m_2,\dots,m_r)\in {\rm Orb}(\s_\gamma/\Z_\ell)$ exists
if and only if $|{\rm Epi}_{\rm o}(\pi_1(\Omega),\Z_\ell)|\neq 0$ (see~\cite{Ha66,MM98}).
Now we can derive an enume\-ra\-ti\-ve counterpart of Theorem~\ref{harv}.

\begin{theorem}\label{van2}
There exists an orbifold\/
$\Omega(g;m_1,m_2,\dots,m_r)\in {\rm Orb}(\s_\gamma/\Z_\ell),$ where
$\ell\geq 2$ and $m_j\geq 2$ for all $j=1,\dots,r$,
if and only if its parameters satisfy the Riemann--Hurwitz equation\ {\rm (RH)}
and none of the following conditions is valid:

\begin{enumerate}
\item[(E1)] $m\nmid\ell$;
\item[(E2)] $g=0$ and $\ell>m$;
\item[(E3)] $s(p)=1$ for some odd prime $p|m$;
\item[(E4)] $2|m$ and $s(2)$ is odd.
\end{enumerate}
\end{theorem}

Proof.
According to~(30), the inequality $|{\rm Epi}_{\rm o}(\pi_1(\Omega),\Z_\ell)|\neq 0$
is equivalent to the condition
$$
 {m^{2g}\phi_{2g}(\ell/m)\cdot E(m_1,m_2,\dots,m_r)}\neq 0.
$$
Now, $m^{2g}\neq 0$ and $\phi_{2g}(d)=0$ if and only if $m\nmid\ell$ or
$g=0$ and $d>1$ (see~(27)). It is clear from~(21) and Corollary~\ref{van1}~(2),
that $E_p(m_1,m_2,\dots,m_r)=0$ if and only if one of
conditions~(E3) and~(E4) of the theorem holds.
\proofend

We can compare Theorem~\ref{van2} with Harvey's Theorem~\ref{harv}.
(E1) and~(E2) are equivalent to~(H2). It is obvious that~(E3) is
equivalent to~(H1) and~(E4) is equivalent to~(H4). As to condition~(H3),
by Corollary~\ref{1-2-3}~(1), $r\neq 1$. Now for $r=2,$ by Corollary~\ref{1-2-3}~(2),
we would have $m_1=m_2=m$ but this contradicts~(RH) for $g=0$ and
$\gamma\neq 1.$ Finally, $r=0$ (in which case $m=1$) is also impossible by~(RH) for
$g=0$ and $\gamma\neq 0.$

\begin{rem}\label{hist} 
{\rm {\bf Additional references.} A related but different enumeration problem with respect
to automorphisms groups (not necessarily cyclic) of orientable surfaces has been considered
and solved by C.\,Maclachlan and A.\,Miller~\cite{MM98}. A simple formula
(with an implicit participation of the function $\phi_k(n)/\phi(n)$) that connects
$|{\rm Epi}(\pi_1(\Omega),\Z_\ell)|$ with the number of homomorphisms
$|{\rm Hom}(\pi_1(\Omega),\Z_\ell)|$ has been obtained by G.\,Jones~\cite{Jo95}.
Diffe\-rent number-theoretic aspects of automorphisms of surfaces partially
related to Harvey's theorem were investigated by W.\,Chrisman~\cite{Ch06};
cf. also the paper~\cite{Si07} by M.\,Sierakowski. In general, the enumerative
approach developed here may have certain useful interactions with combinatorics
of automorphism groups of Riemann surfaces investigated in numerous recent and
older publications (cf., e.g.,~\cite{CI04,Gi76,Wo07}).}
\end{rem}

\subsection{Map enumeration}
A (topological) {\it map}\, is a proper cell embedding of a finite
connected planar graph (generally with loops and multiple edges) in an
orientable surface. A map is called {\it rooted} if an edge-end (a
vertex-edge incidence pair), is distinguished in it as its root. Unlike maps
without a distinguished root (called {\it unrooted\/} maps), rooted maps
do not have non-trivial automorphisms.
Unrooted maps are considered here up to orientation-preserving isomorphism;
i.e., as maps on a surface with a distinguished orientation.
For a general combinatorial
and algebraic theory of maps see~\cite{CoM92,JS78}. Most of numerous results
obtained so far for counting maps are concerned with rooted {\sl planar} maps,
that is, rooted maps on the sphere. Typically maps are counted with respect to
the number of edges.

Let $(g;m_1,m_2,\dots,m_r)$ be a signature of an orbifold and $b_i\geq 0,\ i\geq 2,$
denote the number of branch points with $m_j=i$. That is, up to reordering,
$(g;m_1,m_2,\dots,m_r)
=(g;\underbrace{2,2,\ldots,2}_{b_2},\underbrace{3,\ldots,3}_{b_3},\ldots)=
[g;2^{b_2}3^{b_3}\ldots\ell^{b_\ell}],$ where brackets indicate to the use
of the parameters $b_i$ rather than $m_j.$ In these terms,
the following theorem is valid:
\medskip

\begin{theorem}~{\rm\cite{MeN04}}\label{umaps}
The number of unrooted maps with\/ $n$ edges on a closed
orientable surface\/ $\s_\gamma$ of genus\/ $\gamma$ is
$$
\hspace*{-1ex}\Theta_\gamma(n)\!
=\!\frac{1}{2n}\sum_{\ell|2n}\!
\sum_{\begin{array}{c}\Omega
\end{array}}
|{\rm Epi}_{\rm o}(\pi_1(\Omega),\Z_\ell)|\sum_{s=0}^{b_2} {2n/\ell\choose s}
{\frac{n}{\ell}\!-\!\frac{s}{2}\!+\! 2\!-\!2g\choose b_2\!-\!s,b_3,
\dots,b_\ell}{\cal N}_g\left(\frac{n}{\ell}\!-\!\frac{s}{2}\right),\eqno(36)
$$
where\/ $\ell$ runs over the divisors of\/ $2n,$\,
$\Omega=\Omega[g;2^{b_2}\!\ldots\ell^{b_\ell}]$\/ runs over
the orbifolds in\/ ${\rm Orb}(\s_\gamma/\Z_\ell),$\/
${\cal N}_g(n)$ denotes the number of rooted maps with\/ $n$ edges on
an orientable surface of genus\/ $g$ (with\/ ${\cal N}_g(n):=0$ if\/ $n$
is not integer) and\/ $\displaystyle{n\choose b_2,b_3,\ldots,b_\ell}$ is
the multinomial coefficient.
\end{theorem}

A similar formula holds for non-isomorphic hypermaps~\cite{MeN07,MeN09}.

The orbicyclic function $E(m_1,m_2,\ldots,m_r)$ participates in~(36) through
formula~(30). As for counting arbitrary rooted maps on orientable surfaces we refer
to the papers~\cite{BC91,AG99} (for any genus $g$, there exists a closed formula,
which is remarkably simple, sum-free, in the planar case and becomes more and more
cumbersome as $g$ grows); see also ~\cite{WL72}.

\medskip
Theorem~\ref{umaps} is a far-reaching generalization of the formula derived by the
author for arbitrary unrooted planar maps~\cite{Li81} (see also~\cite{Li98}).
Combinatorially, the proof of~(36) follows the general reductive scheme ela\-bo\-rated
for planar maps: relying upon Burnside's (orbit counting) lemma, the enu\-me\-ra\-tion
of unrooted maps (of a certain class) is reduced to the enumeration of their
corres\-pon\-ding rooted {\it quotient maps} (i.e., orbifold maps) with respect
to all possible orientation-preserving automorphisms
of the underlying surface. In the case of arbitrary planar maps, quotient maps
are merely almost arbitrary planar maps themselves, and we had nothing to do
with the functions $|{\rm Epi}_{\rm o}|$ and $E(m_1,m_2,\ldots,m_r)$ (naturally,
we made use of Euler's $\phi(m)$ in place of them). Burnside's lemma explains
the r\^ole of cyclic groups in Th.~\ref{umaps}.

\medskip
Given $\gamma$, the above formulae make it possible to list all
admissible orbifolds
$\Omega=\Omega(g;m_1,m_2,\ldots,m_r)\in{\rm Orb}(\s_\gamma/\Z_\ell)$.
In view of~(33) and~(34), for $\gamma\geq 2,$ the set
$\bigcup_{\ell\geq 1} {\rm Orb}(\s_\gamma/\Z_\ell)$ is finite.
On the contrary, for $\gamma=0,1,$ there are two infinite families of
or\-bi\-folds, namely, ${\rm Orb}(\s_0/\Z_\ell)=\{\Omega(0;\ell,\ell)\}$
and ${\rm Orb}(\s_1/\Z_\ell)=\{\Omega(1;\emptyset)\},\,\ell=1,2,\ldots$;
(of course, in~(36) they are restricted to the finite set of divisors of $2n$).
Moreover, there are four other orbifolds for $\gamma=1$ with $g=0$ and
$\ell=2,3,4,6.$ Denote
$$
A(\gamma):=\big|\bigcup_{\ell\,\geq 1} {\rm Orb}(\s_\gamma/\Z_\ell)\big|,
$$
where ${\rm Orb}(\s_\gamma/\Z_1)=\{\Omega(\gamma;\emptyset)\},$ and let
$A_g(\gamma)$ denote the cardinality of the subset of all orbifolds
$\Omega(g;m_1,m_2,\ldots,m_r)$ of genus $g$ in this set.
In particular~\cite{Br90,MeN04,Ku90},
$$
A(2)=10,A(3)=17,A(4)=25,A_0(2)=8,A_0(3)=12,A_0(4)=18.
$$
Apparently these functions could be investigated in general basing upon
the restrictions considered above.

\subsection*{Concluding remarks}
{\bf (1)} No publication on multivariate arithmetic functions constructed
similarly to (1) (see formula~(15)) was known to me until recently when
the preprint~\cite{Mi09} appeared, where its author N.\,Minami
considered (in quite a different context) the rational-valued function
$\frac{1}{m}\sum_{k=1}^m(k,m_1)(k,m_2)\cdots(k,m_r)$
and derived a simple formula for its calculation.

\medskip
{\bf (2)} Addressing the polynomial $h_s(x)$ and formula (10)
it would be interesting to find an explanation to the following
observation (belongs to Alexander Mednykh, oral communication):
$h_s(x)=\frac{\chi({\bf C}_s,x)}{x(x-1)}$, where $\chi({\bf C}_s,x)$
denotes the {\sl chromatic polynomial} of a cycle of length $s$
(concerning the chromatic polynomials of graphs we refer the reader,
e.g., to Chapter~IX of the monograph~\cite{Tu84}, in particular, to
the formula for the chromatic polynomials of the cycles given in Theorem IX.24).

\end{document}